\documentclass[11pt, a4paper, reqno]{amsart}

\usepackage{amsfonts,amsmath,amsthm}
\usepackage{amssymb}
\usepackage{latexsym}

\usepackage[english]{babel}

\numberwithin{equation}{section}

\newtheorem{theorem}{Theorem}[section]

\newtheorem{lemma}[theorem]{Lemma}
\newtheorem{corollary}[theorem]{Corollary}

\theoremstyle{definition}

\theoremstyle{remark}

\newtheorem{remark}[theorem]{Remark}

\newcommand{\Id}{\mathrm{Id}}

\newcommand{\dopu}{{:}\allowbreak\ }

\newcommand{\rest}[2]{#1\raisebox{-0.3ex}{\mbox{$\mid_{#2}$}}}

\newcommand{\loglike}[1]{\mathop{\rm #1}\nolimits}

\newcommand{\sign}{\loglike{sign}}

\newcounter{abc}   
\newcounter{iiiii} 

\newenvironment{aequivalenz}
{\setcounter{iiiii}{0}
\begin{list}%
{{\rm (\roman{iiiii})}}
{\usecounter{iiiii}
\parsep=0pt plus 1pt
\topsep=1pt plus 2pt minus 1pt
\itemsep=1pt plus 2pt minus 1pt
\leftmargin=3\baselineskip \labelsep=.6\baselineskip
\labelwidth=2.4\baselineskip
\rightmargin 0pt}%
}
{\end{list}}

\newenvironment{statements}%
{\setcounter{abc}{0}
\begin{list}%
{{\rm (\alph{abc})}}
{\usecounter{abc}
\parsep=0pt plus 1pt
\topsep=1pt plus 2pt minus 1pt
\itemsep=1pt plus 2pt minus 1pt
\leftmargin=3\baselineskip \labelsep=.6\baselineskip
\labelwidth=2.4\baselineskip
\rightmargin 0pt}%
}
{\end{list}}

{\begin{list}%
{$\bullet$}
{\leftmargin=3\baselineskip
\labelsep=\baselineskip \labelwidth=2.5\baselineskip
\parsep=0pt plus 1pt
\topsep=1pt plus 2pt minus 1pt
\itemsep=1pt plus 2pt minus 1pt
\rightmargin 0pt}%
}
{\end{list}}

\newcommand{\bea}{\begin{eqnarray*}}
\newcommand{\eea}{\end{eqnarray*}}
\newcommand{\beq}{\begin{equation}}
\newcommand{\eeq}{\end{equation}}
\newcommand{\begsta}{\begin{statements}}
\def\endsta{\end{statements}}
\newcommand{\begaeq}{\begin{aequivalenz}}
\def\endaeq{\end{aequivalenz}}

\def\DP{Daugavet property}

\begin{document}

\title[The geometry of  $L^p$-spaces]{The geometry of  $L^p$-spaces 
over atomless measure spaces and the Daugavet property}

\author{Enrique A. S\'anchez P\'erez and Dirk Werner}

\date{\today. \textit{File}: \texttt{\jobname.tex}}



\thanks{The first author was partially
 supported by a grant from the Generalitat Valenciana (BEST/2009/108)
 and a grant from the Universidad Polit\'ecnica de Valencia
 (PAID-00-09/2291). Support of the Ministerio de Educaci\'on y 
Ciencia, under project \#MTM2006-11690-C02-01 (Spain) and FEDER is also
gratefully acknowledged.}

\address{Instituto
Universitario de Matem\'atica Pura y Aplicada,  Universidad
\linebreak
Polit\'ecnica de Valencia, Camino de Vera s/n, 46071 Valencia, Spain.}
\curraddr{Department of Mathematics, Freie Universit\"at Berlin,
Arnimallee~6, D-14\,195~Berlin, Germany. }
\email{easancpe@mat.upv.es}

\address{Department of Mathematics, Freie Universit\"at Berlin,
Arnimallee~6, \qquad {}\linebreak D-14\,195~Berlin, Germany}
\email{werner@math.fu-berlin.de}

\begin{abstract}
We show that $L^p$-spaces over atomless measure
spaces can be characterized in terms of a $p$-concavity type geometric property
that is related with the Daugavet property.
\end{abstract}

\maketitle

\thispagestyle{empty}

\section{Introduction}

A Banach space $Y$ is said to have the \textit{Daugavet property} if for every
rank one 
operator
$T\dopu Y \to Y$, the \textit{Daugavet equation}
$$
\|\Id+T\|= 1 + \|T\|
$$ 
is satisfied; in this case, it is known that the
equation is satisfied for every weakly compact operator 
(\cite[Theorem~2.3]{ams2000}). 
Although for $L^1(\mu)$ spaces over an atomless measure $\mu$ 
this property is always fulfilled, it is known that this equation
is 
only satisfied for a compact  operator $T$
on $L^p$ for $1<p<\infty$ when its norm is an eigenvalue of $T$; this
result can be extended to  uniformly convex or uniformly smooth Banach
spaces, and also to locally uniformly convex Banach spaces (see
Corollary~2.4 and Theorem~2.7 in \cite{abraali0} or Section~4 in
\cite{ams2000}). 

After this negative
result, some efforts have been made in order to find a similar lower
estimate for $\| \Id + T\| $ in terms of $\|T\|$ in general Banach
spaces or for the particular case of $L^p$ spaces. Based on the early
ideas of Benyamini and Lin in \cite{benyalin} several authors have
been working in the direction of finding nice lower bounds for
$\|\Id+T\|$ in terms of a function $\psi\dopu (0,+\infty) \to (1,+\infty)$
such that the inequality $\|\Id+T\| \ge 
\psi(\|T\|)$ holds for all compact operators $T\dopu Y \to Y$ (see for instance
\cite{boyK,martin,maroik,oik,popov,schep} and \cite{KMM-Indiana} 
and the references therein).  
As  in these cases, we are interested in this
paper in finding a good alternative to the Daugavet equation for $L^p$
spaces, or in a 
more general sense, for Banach function spaces satisfying certain
$p$-convexity type requirements. Our aim is to give a geometrical
description of $L^p(\mu)$ spaces defined over measures $\mu$ without
atoms in the same geometrical terms as for
spaces with the Daugavet property (slices and the geometry of the
unit ball). In order to do that, we use $p$-convexity and
$p$-concavity properties of quasi-Banach function spaces for
developing a sort of $p$-convexification technique that allows us to
obtain the desired geometrical description. Regarding the Daugavet
property for Banach function spaces the results that are nowadays
known are in a sense negative; for instance, in the class of 
Orlicz spaces over atomless finite measure spaces, the spaces
that satisfy the Daugavet property with respect to the Luxemburg norm are 
isomorphic to $L^1$
(see \cite[Theorem~2.5]{acosta}). However, it must be noted that there
are Banach function spaces other than $L^1(\mu)$ and $L^\infty(\mu)$
over atomless measures $\mu$ that satisfy the Daugavet property (see
for instance Section~5 in 
\cite{bilik}; an explicit example is $c_0(L^1[0,1])$,
the $c_0$-sum of $L^1[0,1]$-spaces).

Let us start by recalling some well-known facts and 
by introducing some notation.
Let $1 \le p < \infty$. A $p$-convex and $p$-concave Banach
lattice can be identified isomorphically and in order with an
$L^p$-space;  if the
corresponding $p$-convexity and $p$-concavity constants are equal to $1$,
then this identification is given by an isometry (see for instance
\cite[Theorem~2.7]{meyer-nieberg}).
In this paper we provide a Daugavet type geometric property which is more
restrictive than the $p$-concavity that is only satisfied for
$L^p$-spaces over measure spaces without atoms (see
Theorem~\ref{carac}). 
In fact, it characterizes this class of spaces.

We remark that we deal with a different $p$-version of the \DP\ in our
paper~\cite{EnrDir1}.

If $Y$ is a Banach space, we
denote as usual by $B_Y$ and $S_Y$ the (closed) unit ball and the
unit sphere respectively. $Y^*$ stands for its dual space.
The slice $S(y^*, \varepsilon)$ of $B_Y$ defined by $y^* \in
B_{Y^*}$ and $\varepsilon >0$ is given by
$$
S(y^*, \varepsilon)= \{y \in B_Y\dopu  \langle y, y^* \rangle \ge 1-
\varepsilon \}.
$$
Notice that for the slice to be non-trivial it is enough to require
that 
$y^* \in
S_{Y^*}$.
Recall that $Y$ has the Daugavet
property if and only if the following geometric property is fulfilled:
for 
every $y \in S_Y$, every $y^* \in S_{Y^*}$ and
every $\varepsilon >0$, there is an element $x \in
S(y^*,\varepsilon)$ such that $\| y+x \| \ge 2- \varepsilon$ (see
\cite[Lemma~2.1]{ams2000}, \cite[Lemma~2.2]{ams2000} or 
\cite[Theorem~2.2]{bulpolis2008}). The reader can find more information on the
geometric description of the Daugavet 
property in \cite{bulpolis2008, ams2000, studia2001} and in 
\cite[Ch.~11]{abraali}.

We also use  standard notation regarding quasi-Banach function spaces. A quasi-Banach space $(E, {\|\cdot\|}_E)$ is a linear space
that is complete with respect to the topology induced by a quasi-norm ${\|\cdot\|}_E$. If $E$ is also a linear lattice,
we say that $(E, {\|\cdot\|}_E)$ is a quasi-Banach lattice if ${\|\cdot\|}_E$ is a lattice quasi-norm
in $E$, i.e., $\| x\|_E \le \|y\|_E$ whenever $x,y \in E$ and $|x| \le |y|$.
Let $(\Omega, \Sigma, \mu)$ be a measure space. A quasi-Banach function
space $X(\mu)$ over the measure $\mu$  is an ideal of $L^0(\mu)$, the
usual $\mu$-a.e.\ 
order is considered, which  is a  quasi-Banach space with a
lattice quasi-norm ${\|\cdot\|}$ such that for every $A \in  \Sigma$ of
finite measure, $\chi_A \in X(\mu)$ (see for instance \cite[Chapter~1.b]{lint}, \cite[Chapter~2.6]{meyer-nieberg} and
\cite[Chapter~2]{libro}
  for definitions and main results regarding these structures, but notice that the last property is not required in some of
 these references). In the case that ${\|\cdot\|}$ is a norm, we say
 that $(X(\mu),{\|\cdot\|})$ is a Banach function space, 
see \cite[Definition~1.b.17]{lint}. We shall write simply $X$ for $X(\mu)$ if the
 measure is fixed in the context.

Let $0<p < \infty$. A quasi-Banach lattice $E$ is \textit{$p$-convex} 
if there is a constant $K$ such
that for every finite sequence $(x_i)_{i=1}^n$ in $E$,
$$
\Bigl\| \Bigl( \sum_{i=1 }^n |x_i|^p \Bigr)^{1/p} \Bigr\|_E
 \le K \Big(\sum_{i =1 }^n \|x_i\|_E^p \Big)^{1/p}.
$$
A quasi-Banach lattice $E$ is \textit{$p$-concave} 
if there is a constant $k$
such that for every finite sequence $(x_i)_{i=1}^n$ in $E$,
$$
\Bigl(\sum_{i =1}^n \|x_i\|_E^p \Bigr)^{1/p} \le  k \Bigl\| \Bigl(
\sum_{i = 1}^n |x_i|^p \Bigr)^{1/p} \Bigr\|_E.
$$
The best constants in
the inequalities above are denoted by $M^{(p)}(E)$ and $M_{(p)}(E)$, 
respectively, and are called the $p$-convexity and the $p$-concavity constants of $E$.
Throughout the paper we will assume for a $p$-convex
Banach function space   that in fact $M^{(p)}(E)=1$, and when
$p$-concavity is required 
that $M_{(p)}(E)=1$. To indicate this, we will say that they are constant~$1$ $p$-convex or constant~$1$ $p$-concave, respectively.

Let $0< p < \infty$. Consider  a quasi-Banach function space
$X(\mu)$. Let
$$
X(\mu)_{[p]}:=\{h \in L^0(\mu)\dopu  |h|^{1/p} \in X(\mu) \}
$$
be its $p$-th power, which is a quasi-Banach function space when
endowed with the quasinorm $\|h\|_{X_{[p]}}:=\| |h|^{1/p} \|^p_X$, $h
\in X_{[p]}$ (see \cite[Ch.~2]{libro}). For $p \ge1$, if $X$ is $p$-convex and $M^{(p)}(X)=1$, then $X(\mu)_{[p]}$ is a Banach
function space, since in this case ${\|\cdot\|}_{X_{[p]}}$ is a norm. If the Banach function space  is $p$-convex, but
the $p$-convexity constant is not~$1$, then ${\|\cdot\|}_{X_{[p]}}$ is
not a norm, but it is equivalent to a norm (see for instance
\cite[Prop.~2.23]{libro}). It is also well known that every $p$-convex
Banach lattice can be renormed in such a way that the new norm is a
lattice norm with $p$-convexity constant equal 
to~$1$ (see \cite[Proposition~1.d.8]{lint}).

\section{Banach function spaces with $p$-th powers having the 
Daugavet property}

Let $0 < p < \infty$ and let $X(\mu)$ be a constant~$1$ $p$-convex
quasi-Banach function space. 
If $f \in X(\mu)$, we can always write it as $f = \sign \{f\}
|f|$. This allows us to define the  (obviously non-linear)
map $i_p\dopu X(\mu) \to X(\mu)_{[p]}$ by 
$$
i_p(f)= \sign \{f\} |f|^p.
$$
Throughout the paper, we shall write 
$$
f^p:=i_p(f)
$$ 
for the sake of simplicity, but notice that for even integers 
$f^p$ is  in general not $|f|^p$. The map $i_p$  is bijective and satisfies 
\begin{equation}\label{eq1}
\|i_p(f)\|_{X_{[p]}} = \| {\sign \{f\}} |f|^p \|_{X_{[p]}} = \|f\|_X^p, \quad f \in X(\mu).
\end{equation}

The inverse map $i^{-1}_{p}\dopu  X_{[p]} \to X$ coincides with $i_{1/p}\dopu  Y \to Y_{[1/p]}$, where $Y=X_{[p]}$.

In what follows we characterize the $p$-convex
Banach function spaces whose $p$-th powers satisfy the Daugavet
property. The key idea to achieve this is to
introduce the notions
 of $1/p$-th
power of a slice and $p$-convexification of an operator
$T\dopu X_{[p]} \to X_{[p]}$.

If  $X(\mu)$ is a constant~$1$ $p$-convex Banach function space, let
$S_{[p]}(x^*,\varepsilon)$ be a slice in $X(\mu)_{[p]}$, where
$x^* \in B_{(X(\mu)_{[p]})^*}$. Consider the  set
$$
S_{[p]}^{1/p}(x^*,\varepsilon):= \{ f \in X(\mu)\dopu   f^p \in S_{[p]}(x^*,\varepsilon) \}.
$$
We call it the \textit{$1/p$-th power of the
slice} $S_{[p]}(x^*,\varepsilon)$.

If $T\dopu X_{[p]} \to X_{[p]}$ is an operator, we define its
\textit{$p$-convexification} $\varphi_T\dopu X(\mu) \to X(\mu)$ by
$$
\varphi_T(f):= i^{-1}_{p} \circ T \circ i_{p}(f) =
(T(f^p))^{1/p}, \quad f \in X.
$$
We also define  $\| \varphi_T\|:= \sup_{f \in B_{X}}
\|\varphi_T(f)\|$. Notice that
\bea
\| \varphi_T\| &=& \sup_{f \in B_{X}} \|\varphi_T(f)\|_X
= \sup_{f \in B_{X}} \|(T(f^p))^{1/p}\|_X  \\
&=&
\sup_{h \in
B_{X_{[p]}}} \|T(h)\|^{1/p}_{X_{[p]}} = \|T\|^{1/p}.
\eea

The following two lemmas provide a geometric description 
of the Daugavet property for a Banach function space
in terms of slices of  the
$p$-convexification.   Their proofs follow the lines of the
ones of Lemmas~2.1, 2.2 and 2.8 in \cite{ams2000}. However, we 
spell out  the arguments that prove the
main equivalences with some detail in order to show the role played by the $p$-convexity of the norm of $X(\mu)$.

\begin{lemma} \label{LemmapthDaugavet}
 Let $X(\mu)$ be a quasi-Banach function space  and let $0<{p<\infty}$. 
The following assertions are equivalent:
\begin{itemize}
\item[(0)] The space satisfies the following.
\begin{itemize}
\item[(i)] $(X_{[p]})^* \ne \{0 \}$.
\item[(ii)] For every finite family of rank-one continuous 
operators $T_i\dopu X_{[p]} \to X_{[p]}$,
$$
\sup_{f \in B_X} \Bigl\| \Bigl(\sum_{i=1}^n \varphi_{T_i}(f)^p
\Bigr)^{1/p} \Bigr\| 
\le \Bigl(\sum_{i=1}^n \|\varphi_{T_i}\|^p \Bigr)^{1/p} .
$$
\item[(iii)]
For every rank-one operator $T\dopu X_{[p]} \to X_{[p]}$,
 $$
 \sup_{f \in B_X} \||f^p+ \varphi_T(f)^p |^{1/p} \|_X= (1 + \|\varphi_T\|^p)^{1/p}.
 $$
\end{itemize}
 \item[(1)] $X$ is constant~$1$ $p$-convex and for every rank-one
operator $T\dopu X_{[p]} \to X_{[p]}$,
 \begin{equation}\label{eq00}
 \sup_{f \in B_X} \||f^p+ \varphi_T(f)^p |^{1/p} \|_X= (1 + \|\varphi_T\|^p)^{1/p}.
 \end{equation}

\item[(2)] $X$ is constant~$1$ $p$-convex and for every $f \in S_{X}$, every $x^* \in S_{(X_{[p]})^*}$ and every $\varepsilon >0$
there is an element $g\in S^{1/p}_{[p]}(x^*,\varepsilon)$ such
that
$$
\| |f^p + g^p |^{1/p}\|_X^p \ge 2- 2\varepsilon.
$$
\item[(3)]
$X(\mu)_{[p]}$ if a Banach function space  over $\mu$ with the Daugavet property.
\item[(4)]
$X$ is constant~$1$ $p$-convex and
for every $f \in S_X$ and every slice $S(x^*_0,
\varepsilon_0)$ of $B_{X_{[p]}}$ there is another non-trivial slice $S_{[p]}(x^*_1,
\varepsilon_1) \subset S_{[p]}(x^*_0, \varepsilon_0)$ such that for
every $g \in S_{[p]}^{1/p}(x^*_1,\varepsilon_1)$ the inequality
$$
\| (f^p + g^p)^{1/p} \|^p_{X} \ge 2- \varepsilon_0
$$
holds.
\end{itemize}
\end{lemma}

\begin{proof}
Let us start with (0) $ \Rightarrow $ (1). Take a finite set
$f_1,\dots ,f_n \in X$. Since the dual space $(X_{[p]})^*$ contains a
non-trivial element $x^*$, we can assume that  $x^* \in
S_{(X_{[p]})^*}$, and we can consider the operators $T_i\dopu X_{[p]}
\to X_{[p]}$ given by $T_i(h):= \langle h, x^* \rangle |f_i|^p$. They
are obviously continuous and $\| T_i\|=\|f_i^p\|_{X_{[p]}}= \|
f_i\|_X^p$. Thus, by~(ii), we have
\bea
\Bigl(\sum_{i=1}^n \|f_i\|_X^p \Bigr)^{1/p} 
&=& 
\Bigl(\sum_{i=1}^n \|\varphi_{T_i}\|^p \Bigr)^{1/p}\\
&\ge&
\sup_{f \in B_X} \Bigl\| \Bigl(\sum_{i=1}^n
\varphi^p_{T_i}(f)\Bigr)^{1/p} 
\Bigr\| \\
&=& 
\sup_{f \in B_X} \Bigl\| \Bigl(\sum_{i=1}^n T_i(f^p) \Bigr)^{1/p} \Bigr\|
\\
&=& 
\sup_{f \in B_X} \Bigl\| \Bigl(\sum_{i=1}^n \langle f^p, x^* \rangle
|f_i|^p \Bigr )^{1/p} \Bigr\| \\
&=&
\sup_{h \in B_{X_{[p]}}} ( \langle h, x^* \rangle)^{1/p} 
\Bigl\| \Bigl(\sum_{i=1}^n |f_i|^p \Bigr)^{1/p} \Bigr\|_X.
\eea
Consequently, $X$ is $p$-convex and $M^{(p)}(X)=1$, and so (1) is
obtained. 

For the converse, since $X$ is $p$-convex and has $p$-convexity
constant equal to~$1$, $X_{[p]}$ is a Banach function space (see for
instance \cite[Proposition~2.23]{libro}), and so its dual space is
non-trivial. 
It only remains to prove~(ii). Take any finite set of rank-one operators
$T_i\dopu X_{[p]} \to X_{[p]}$, $i=1,\dots ,n$. Each of them can be
written as $T_i=x_i^* \otimes f_i^p$, where $\|x_i^*\|=1$ and $f_i \in
X$. Then for every $f \in B_X$,
\bea
\Bigl\| \Bigl(\sum_{i=1}^n  (\varphi_{T_i}(f))^p \Bigr)^{1/p} \Bigr\| 
&\le&
 \Bigl\| \Bigl(\sum_{i=1}^n |\langle f, x_i^* \rangle|^p |f_i|^p
 \Bigr)^{1/p} \Bigr\| \\
&\le&
\Bigl(\sum_{i=1}^n \|f_i\|^p \Bigr)^{1/p} = \Bigl(\sum_{i=1}^n \|\varphi_{T_i}\|^p\Bigr)^{1/p}.
\eea
This gives (0).

 Let us now prove the equivalence of (1) and (2). First notice that
(1) is equivalent to the fact that for every rank-one operator
$T\dopu X_{[p]} \to X_{[p]}$,
 $$
 \sup_{g \in B_X} \||g^p+ T(g^p) |^{1/p} \|^p_X=1 + \|T\|.
 $$
 For (1) $\Rightarrow$ (2), take $f \in S_{X(\mu)}$, $x^* \in
S_{(X(\mu)_{[p]})^*}$ and $\varepsilon
>0$. Consider $T=x^* \otimes f^p$.
The equality above can be written as
\begin{equation}\label{eq2}
 \sup_{g^p \in B_{X_{[p]}}} \|g^p+ T(g^p) \|_{X_{[p]}}=1 + \|T\|.
\end{equation} 
In particular, this implies that we can assume by Lemma~11.4
in \cite{abraali} (or \cite[p.~78]{Dirk-IrBull})
that $T$ and hence $x^*$ and $f$ are of norm one. Take $h \in S_X$ such that
$$
\||h^p+ T(h^p) |^{1/p} \|^p_X\ge 2-\varepsilon.
$$
We can also assume that $\langle h^p,x^* \rangle \ge 0$ (otherwise
replace $h$ by $-h$). Notice first that since $X(\mu)$ is
$p$-convex with  constant~$1$, 
$$
1+ \langle h^p, x^* \rangle= \| |h^p|^{1/p} \|^p_X + \|
|T(h^p)|^{1/p}\|^p_X \ge 2-\varepsilon,
$$
which implies that $\langle h^p, x^* \rangle \ge 1-\varepsilon$.
Consequently, $h \in S^{1/p}_{[p]}(g^*,\varepsilon)$. On the other
hand, using again the constant~$1$ $p$-convexity of $X(\mu)$,
\bea
2-\varepsilon 
&\le& \| |h^p+T(h^p)|^{1/p}\|^p_X \\
&\le& \| |h^p+f^p|^{1/p} \|^p_X + \| |T(h^p)-f^p|^{1/p} \|^p_X \\
&= & \| | h^p+f^p|^{1/p} \|^p_X + \| |(\langle h^p,x^*
\rangle-1)f^p|^{1/p} \|_X^p \\
&\le& \| | h^p+f^p|^{1/p} \|^p_X + (1-\langle h^p,x^* \rangle) \\
&\le&
\| | h^p+f^p|^{1/p} \|^p_X +  \varepsilon.
\eea
This gives (2). 

For the converse, we can suppose again that $T$ is
defined as $T= x^* \otimes f^p$ for two norm one elements $x^*$
and $f$. Let $\varepsilon >0$ and $h \in
S_{[p]}^{1/p}({x^*},\varepsilon)$ such that
$$
\| |f^p + h^p |^{1/p}\|_X^p \ge 2- 2\varepsilon.
$$
Then,  by the constant~$1$ $p$-convexity of $X(\mu)$,
\bea
2- 2\varepsilon &\le& \| |f^p + h^p |^{1/p}\|_X^p \\
&=& 
\| |f^p-T(h^p) + T(h^p) + h^p |^{1/p}\|_X^p\\
&\le& 
\| |f^p-T(h^p)|^{1/p}\|_X^p + \||T(h^p) + h^p |^{1/p}\|_X^p \\
&\le& (1-\langle h^p, x^* \rangle) + \||T(h^p) + h^p |^{1/p}\|_X^p \\
&\le& \varepsilon  + \||T(h^p) + h^p |^{1/p}\|_X^p.
\eea
Since this holds for every $\varepsilon >0$ and the
 converse inequality always holds, we obtain the result.

Taking into account the definition of the isometric map $i_{p}$, the
definition of
the norm ${\|\cdot\|}_{X_{[p]}}$ and (\ref{eq2}), 
the equivalence of (1) and (3) becomes  obvious.
Notice that the fact that $X(\mu)_{[p]}$ is a Banach function space  over $\mu$ is equivalent to the fact that
$X(\mu)=(X(\mu)_{[p]})_{[1/p]}$ is constant~$1$ $p$-convex (see for instance \cite[Proposition~2.23(ii)]{libro}).

Similar arguments prove (3) $\Rightarrow$ (4); a direct proof can be given using Lemma~2.1(a) in
\cite{ams2000}, the definition of the norm in $X_{[p]}$
and the fact that every element $h \in X_{[p]}$ can be written as
$f^p$ for some $f \in X$. (4) $\Rightarrow$ (2) is obvious.
\end{proof}

\begin{remark}
$L^p(\mu)$ spaces over a non-atomic measure $\mu$  satisfy the
equivalent statements of Lemma~\ref{LemmapthDaugavet}; this is a
direct consequence of  $(L^p(\mu))_{[p]}=L^1(\mu)$ and the well-known
fact that $L^1(\mu)$ satisfies the Daugavet property (see for instance
\cite[Theorem~3.2]{abraali0}, or the example after Theorem~2.3 in
\cite{ams2000} for a simple proof).  However, we can easily construct 
Banach function spaces which are not $L^p$ spaces but their $p$-th
powers have the Daugavet property. For instance, consider a
$\sigma$-finite  atomless measure space  $(\Omega, \Sigma,\mu)$
and an infinite measurable partition $\{A_i\}$ of $\Omega$ and take a
Banach space $F$ with a $1$-unconditional normalized Schauder basis
endowed with its natural Banach function space structure given by the
pointwise order. Consider the Banach space $X$ defined as the $F$-sum
of the spaces $L^1(\rest{\mu}{A_i})$, where $\rest{\mu}{A_i}$ denotes the
restriction of $\mu$ to $A_i$, $i \in \mathbb{N}$, that is, $X$ is the space
of sequences $(f_i)$ such that  $f_i \in L^1(\rest{\mu}{A_i})$ and $(\|f_i\|)
\in F$. If $F$ has the positive Daugavet property (i.e., every positive
rank one operator on $F$ satisfies the Daugavet equation), then
Theorem~5.1 in \cite{bilik} ensures that the $F$-sum has the Daugavet
property. The spaces $\ell^1$ and $\ell^\infty$ satisfy the positive
Daugavet property, but the reader can find other examples in
\cite[Section 5]{bilik}. It is easy to see that the $1/p$-th power of
$X$ is also a Banach function space  and it can be identified isometrically with the
$F_{[1/p]}$-sum of the spaces $L^p(\rest{\mu}{A_i})$. Since
$(X_{[1/p]})_{[p]}=X$ has the Daugavet property, $X_{[1/p]}$ satisfies 
the assertions of Lemma~\ref{LemmapthDaugavet}.

Other examples can be constructed using the fact that spaces of
Bochner integrable functions over atomless measures  satisfy the
Daugavet property (see again the example after Theorem~2.3 in
\cite{ams2000}). Let $Y$ be a Banach lattice, let $\mu$ be a measure
without atoms and consider the Bochner space $L^1(\mu,Y)$. It is a
Banach lattice when the natural order is considered; assume that it is
also an order continuous Banach lattice with a weak unit. Then it can
be represented as a Banach function space $Z$ (see for instance
\cite[Theorem~1.b.14]{lint}). Since the Daugavet property is preserved
under isometries, we obtain that $Z_{[1/p]}$ satisfies  the
statements of Lemma~\ref{LemmapthDaugavet}.

\end{remark}

\begin{remark}
Note that although the assertions in Lemma~\ref{LemmapthDaugavet} have
been stated in terms of rank one operators, the equivalences also hold
when other classes of operators satisfying the Daugavet equation in
$X_{[p]}$ are considered. Therefore, it includes for instance the
weakly compact operators and further classes, see for example
\cite{ams2000} and \cite{studia2001}.
\end{remark}

The following ``sign  independent'' inequality is crucial for the
computations regarding the $p$-convexification of the Daugavet
property.

Given $1\le p<\infty$, we denote by $p'$ the conjugate exponent
defined by $1/p + 1/p'=1$. Also, we let $k(p)=1$ if $p\ge p'$ and
$k(p)= 2^{(p'/p) - 1}$ if $p<p'$. It follows
$$
(a^{p/p'} + b^{p/p'})^{p'/p} \le k(p) (a+b )
$$
for real numbers $a,b\ge0$.

\begin{lemma} \label{desig}
Let $1 \le p < \infty$ and consider two elements $f$ and $g$ in the
unit ball of the constant~$1$ $p$-convex Banach function space~$X$. Then 
$$
\||f^p-g^p|^{1/p}\|^p \le  \|f-g\|^p+ p (2 k(p))^{p/p'}) \|f-g\|.
$$
Consequently, the map $i_p\dopu X\to X_{[p]}$ mapping $f$ to $f^p$ is
continuous. 
\end{lemma}

\begin{proof}
Let $1 \le p < \infty$ and consider $c^p:= \sign \{c\}|c|^p$ for
 every $c \in \mathbb{R}$. Let $a, b \in \mathbb{R}$. Then
we have to take into account two cases:

1) $\sign \{a\} \ne \sign  \{b\}$. Suppose without loss of generality
that ${a \ge 0}$ and ${b \le 0}$. Then
$$
|a^p-b^p|^{1/p} = |a^p + |b|^p|^{1/p} \le |a + |b||= |a-b|.
$$

2) $\sign \{a\} = \sign  \{b\}$. Then it is known that 
$$
| a^p -b^p|^{1/p} \le \big( p\, |a^{p-1}+b^{p-1}|\cdot |a-b| \big)^{1/p}
$$
(see for instance \cite[Section 2.2]{libro}). 

Take now two functions $f,g \in B_X$ and 
put $A=\{\omega\dopu \sign\{f(\omega)\} \neq \sign\{g(\omega)\}\}$ and 
$B=\{\omega\dopu \sign\{f(\omega)\} = \sign\{g(\omega)\}\}$.
 Then by case~1)
$$
\||f^p-g^p|^{1/p} \chi_A\|^p \le \||f-g| \chi_A\|^p.
$$
Since $p-1=p/p'$, by the H\"older inequality for the Banach lattice
$X$ (see for instance Proposition 1.d.2 in \cite{lint}), we obtain
(see also \cite[Section 2.2]{libro} for the pointwise inequalities
involved) 
\bea
\||f^p-g^p|^{1/p} \chi_B\|^p 
&\le&  
p \, \| \big( |f^{p-1}+g^{p-1}|\cdot |f-g| \big)^{1/p}\chi_B  \|^p \\ 
&\le& 
p \, \| (|f^{p/p'}+g^{p/p'}|^{p'/p})^{1/p'} \cdot |f-g|^{1/p} \chi_B\|^p\\
&\le& 
p \, \| (|f^{p/p'}+g^{p/p'}|^{p'/p}) \chi_B \|^{p/p'} \cdot \||f-g| \chi_B \|\\
&\le&  
p k(p)^{p/p'} \| |f+g| \chi_B \|^{p/p'} \cdot \||f-g| \chi_B \|\\
&\le& 
p (2 k(p))^{p/p'}  \||f-g| \chi_B \|.
\eea
Therefore, since by the constant $1$ $p$-convexity of $X$ the inequality 
$$
\||f^p-g^p|^{1/p}\|^p  \le \||f^p-g^p|^{1/p} \chi_A\|^p + \||f^p-g^p|^{1/p} \chi_B\|^p
$$
is satisfied, we obtain the result.
\end{proof}

\begin{remark}
We can relate our
property for rank one operators with the general $\psi$-Daugavet
property  for Banach spaces that has been quoted in the Introduction
(see \cite{boyK,oik,schep}). For example in Theorem~2.1 of \cite{oik}
inequalities like $\| \Id + T\| \ge (1+ c_p\|T\|^p)^{1/p}$ for a compact
operator $T\dopu X \to X$ are considered, where $c_p$ is a
non-negative constant. In our case we obtain the following similar
estimate in terms of the $p$-convexification $\varphi_T$. For
instance, if $X_{[p]}$ is a Banach function space   with the Daugavet
property and $T\dopu X_{[p]} \to X_{[p]}$ is weakly compact, we obtain 
$$
 (1 + \|\varphi_T\|^p)^{1/p}\le \sup_{f \in B_X} \||f^p+ \varphi_T(f)^p |^{1/p} \|_X 
$$
or equivalently
$$
 (1 + \|T\|)^{1/p}\le \sup_{f \in B_X}\||f^p+ T(f^p)|^{1/p}\|_X.
$$
Clearly in the case of positive operators, and using the estimate given in the proof of Lemma~\ref{desig} for the case of different signs, this inequality gives also
$$
(1 + \|\varphi_T\|^p)^{1/p} \le \sup_{f \in B_X}\|f+ \varphi_T(f)\|_X.
$$
\end{remark}

The following lemma is similar to Lemma~2.8 in \cite{ams2000}.

\begin{lemma} \label{finite}
Suppose that $X_{[p]}$ is a Banach space with  the Daugavet property. 
Then for every
finite dimensional subspace $X_0$ of $X$, every $\varepsilon >0$
and every $x^* \in (X_{[p]})^*$ there is an element $g \in
S^{1/p}_{[p]}(x^*,\varepsilon)$ such that for every $f \in X_0$
and $t \in \mathbb{R}$
$$
\|( (tg)^p + f^p)^{1/p}\|^p \ge (1-\varepsilon)(|t|^p + \|f\|^p).
$$
\end{lemma}

\begin{proof}
Take $\delta>0$ such that $\delta^p+ p(2k(p))^{p/p'} \delta \le
\varepsilon/2$, where $k(p)$ is defined as above, a
finite dimensional 
subspace $X_0$ of $X$ and a finite $\delta$-net $\{f_1,\dots ,f_n\}$
in $S_{X_0}$. Applying Lemma~\ref{LemmapthDaugavet}(4) we find a
sequence of slices $S_{[p]}(x_n^*, \varepsilon_n) \subset \ldots\subset
S_{[p]}(x_1^*, \varepsilon_1)   \subset S_{[p]}(x^*,
\varepsilon)$ such that
$$
\| (f^p_k + g^p)^{1/p} \|^p_X \ge 2- \delta^p
$$
for all $g \in S^{1/p}_{[p]}(x^*_k,\varepsilon_k)$, $k=1,\dots ,n$. If we consider
elements $g$ in $S^{1/p}_{[p]}(x_n^*, \varepsilon_n)$, these
inequalities are true for all $k=1,\dots ,n$. Consequently, by 
Lemma~\ref{desig} and the constant~$1$ $p$-convexity of $X$, for every
$g \in S^{1/p}_{[p]}(x_n^*, \varepsilon_n)$ and $f \in S_{X_0}$ 
there is an index $k \in \{1,\dots ,n\}$ such that
\bea
\| (f^p+g^p)^{1/p} \|^p_X 
&\ge& \| (f_k^p+g^p)^{1/p} \|^p_X - \| |f_k^p-f^p|^{1/p} \|^p_X \\
&\ge&  2- \delta^p - \varepsilon/2 \ge 2- \varepsilon.
\eea

Now, if $0 \le s \le t$ are real numbers such that $t^p+s^p=1$, then
for all such $g$ and $f$,
\bea
\| ((t g)^p + (s f)^p)^{1/p} \|^p 
&=& \| (t^p(g^p+f^p) - |s^p-t^p|f^p)^{1/p} \|^p \\
&\ge&  t^p\|(g^p+f^p)^{1/p}\|^p - |s^p-t^p|\|f \|^p \\
&\ge& t^p(2-\varepsilon) +s^p-t^p = 1- \varepsilon.
\eea
Since the same calculations can be done for $t \le s$, we obtain
the following inequality for every $t \ge0$ and $f \in X_0$:
\bea
\| ((t g)^p + f^p)^{1/p} \|^p
&=& \Bigl\| \Bigl((t g)^p + \Bigl(\|f\|
\frac{f}{\|f\|}\Bigr)^p\Bigr)^{1/p} \Bigr\|^p  \\
&=& \Bigl\| \Bigl(\frac{t^p g^p}{t^p+ \|f\|^p} +
\frac{\|f\|^p(f/\|f\|)^p}{t^p+ \|f\|^p} \Big)^{1/p} \Bigr\|^p \cdot
(t^p + \|f\|^p ) \\
&\ge&  (1-\varepsilon)(t^p + \|f\|^p ).
\eea
The symmetry of the norm allows to obtain the same inequality for
every $t \in \mathbb{R}$, replacing $t$ by $|t|$.
\end{proof}

The following lemma makes it clear that the Daugavet type equation
(\ref{eq00}) fails in the presence of atoms. 

\begin{lemma} \label{atom}
Let $(\Omega,\Sigma,\mu)$ be a measure space. Let $0<p<\infty$,
let $X(\mu)$ be a constant~$1$ $p$-convex quasi-Banach
function space over $\mu$ and suppose that $\mu$ has an atom. Then
there is a rank one operator $T\dopu X_{[p]} \to X_{[p]}$ such that 
$$
\sup_{f \in B_X} \| (f^p + T(f^p))^{1/p} \|_X < (1+\|T\|)^{1/p}.
$$
\end{lemma}

\begin{proof}
Recall that by the constant~$1$ $p$-convexity of $X$, $X_{[p]}$ is
a Banach function space.
 Let $\{a\}$ be an atom for $\mu$. Then $0 < \mu(\{a\}) < \infty$ and the characteristic function
$\chi_{\{a\}}$ belongs  to $X_{[p]}$, and defines a (continuous)
functional of $(X_{[p]})^*$ by $\langle h, \chi_{\{a\}}\rangle =
 \int \chi_{\{a\}} h \,d\mu= h(a) \mu(\{a\})$, $h \in X_{[p]}$.
Let $T$ be the non-null rank one operator $T=- \frac{\chi_{\{a\}}}{\mu(\{a\})} \otimes \chi_{\{a\}}$. Then
\bea
\sup_{f \in B_X} \| ( f^p +T(f^p))^{1/p} \| 
&=& \sup_{f \in B_X} \|
( f^p(a)\chi_{\{a\}}  + f^p\chi_{\Omega \setminus \{a\}} -
f^p(a)\chi_{\{a\}})^{1/p} \| \\
&=& \sup_{f \in B_X} \| f \chi_{\Omega \setminus \{a\}} \| \le 1 < (1+ \|T\|)^{1/p},
\eea
as claimed.
\end{proof}

The following result provides the desired geometric characterization of $L^p$ spaces over atomless
measure spaces. Recall that an abstract $L^p$ space is a Banach lattice $E$ for which
for every couple of disjoint elements $x,y \in E$, the equality
$\|x+y\|^p = \|x\|^p + \|y\|^p$
holds (see for instance \cite[Definition~1.b.1]{lint}).

\begin{theorem} \label{carac}
Let $1 \le p < \infty$ and let $X(\mu)$ be a quasi-Banach function
space over $\mu$. The following statements are equivalent:

\begin{itemize}
 \item[(i)] $X$ is an abstract $L^p$ space such that $X_{[p]}$ has the Daugavet property.

\item[(ii)] $X$ is equal to $L^p(h \,d\mu)$  for some $0<h \in
  L^1(\mu)$ and the measure $\mu$ does not have atoms.

\item[(iii)] $(X_{[p]})^* \ne \{0\}$, $X$ is constant~$1$ $p$-concave,
  and for every finite set of operators $T_i\dopu X_{[p]} \to
  X_{[p]}$, $i \in \{1,\dots ,n\}$ we have that 
$$
\sup_{f \in B_X} \Bigl\| \Bigl(\sum_{i=1}^n \varphi_{T_i}(f)^p
\Bigr)^{1/p} \Bigr\| 
\le \Bigl(\sum_{i=1}^n \|\varphi_{T_i}\|^p \Bigr)^{1/p} 
$$
and
$$
 \sup_{f \in B_X} \||f^p+ \varphi_{T_i}(f)^p |^{1/p} \|_X= (1 + \|\varphi_{T_i}\|^p)^{1/p}.
$$

\item[(iv)] $X$ is constant~$1$ $p$-convex, constant~$1$ $p$-concave and for every rank one operator
$T\dopu X_{[p]} \to X_{[p]}$,
$$
\sup_{f \in B_X} \||f^p+ \varphi_{T}(f)^p |^{1/p} \|_X= (1 + \|\varphi_{T}\|^p)^{1/p}.
$$

 \item[(v)]
$X$ is constant~$1$ $p$-convex and
 for every slice $S_{[p]}(x^*,\delta)$, every $\varepsilon >0$ and every finite dimensional subspace $X_0$ of $X$
there is an element $g \in S^{1/p}_{[p]}(x^*,\delta)$ such that
for every $f_1,\dots ,f_n \in X_0$ and $\alpha_i \ge0$ satisfying
$\sum_{i=1}^n \alpha^p_i=1$,
$$
\Bigl\|\Bigl( \sum_{i=1}^n |(\alpha_i g)^p + f_i^p|\Bigr)^{1/p}
\Bigr\|^p 
\ge (1- \varepsilon)\Bigl(  \|g\|^p  + \sum_{i=1}^n\|f_i\|^p\Bigr).
$$

\item[(vi)] $X$ is constant~$1$ $p$-convex and for every slice $S_{[p]}(x^*,\delta)$, every $\varepsilon >0$ and  every finite dimensional subspace $X_0$ of $X$
there is an element $g \in S^{1/p}_{[p]}(x^*,\delta)$ and an
element $x_0^* \in B_{(X_{[p]})^*}$ such that
$$
\Bigl(\frac{|f^p+g^p|}{\|f\|^p + \|g\|^p} \Bigr)^{1/p}  
\in S^{1/p}_{[p]}(x_0^*,\varepsilon)
$$ 
for every $f \in X_0$.
\end{itemize}
\end{theorem}

\begin{proof}
(i) $\Rightarrow$ (ii). Since $X$ is an abstract $L^p$ space and $1 \le p
< \infty$, $X$ is  in particular a $\sigma$-order continuous Banach function space. Using
 a Maurey-Rosenthal type factorization argument (see for example
Corollary~6.17 in \cite[Ch. 6]{libro}), we find that there is a
function $0<g$ such that the identity map $\Id\dopu X(\mu) \to X(\mu)$
factorizes through $L^p(\mu)$ by means of the multiplication
operators $M_g\dopu X \to L^p(\mu)$ and $M_{1/g}\dopu L^p(\mu) \to X$; in
fact, $M_g$ defines an isometry (notice that for applying the
Corollary 6.17 quoted above it is necessary to take into account that
the operator $M_g$ always has dense range). Therefore, the space
$X(\mu)$ can be identified isometrically and in order with $L^p(h\,
d\mu)$, $h=g^p$, and its elements are the same functions. So,
$X_{[p]}=L^1(h \,d\mu)$ has the Daugavet property, and therefore by
Lemma~\ref{atom} with $p=1$ the measure $h \,d\mu$ does not have  atoms.
Consequently, $\mu$ does not have  atoms either.

For (ii) $\Rightarrow$ (i), just recall that an $L^1(\nu)$-space
over an atomless measure space has the Daugavet property (see
for instance \cite[Theorem~3.2]{abraali0} or the example after
Theorem~2.3 
in \cite{ams2000}).

By Lemma~\ref{LemmapthDaugavet}, (ii) implies the equivalent
statements (iii) and (iv), taking into account that $X$ can be
written as a Banach function space  over the measure $h \,d\mu$. Clearly, (iv)
implies (i).

Let us now show (i) $\Rightarrow$ (v).  Assume that $X$ is an
$L^p$-space  and $X_{[p]}$ has the Daugavet property. 
Then 
Lemma~\ref{finite} provides  for every finite dimensional subspace
$X_0$ of $X$, every $\varepsilon
>0$ and every $x^* \in (X_{[p]})^*$  an element $g \in
S^{1/p}_{[p]}(x^*,\varepsilon)$ such that for every $f \in X_0$
and $t \in \mathbb{R}$
$$
\|( (tg)^p + f^p)^{1/p}\|_X^p \ge (1-\varepsilon)(|t|^p +
\|f\|_X^p).
$$

Thus, taking into account that $X$ is an $L^p$-space (and then
also constant~$1$ $p$-concave), for every finite set of elements
$f_1,\dots ,f_n \in X_0$ and positive real numbers $\alpha_i$
such that $\sum_{i=1}^n \alpha^p_i=1$, we obtain
\begin{eqnarray}
\Bigl\| \Bigl(\sum_{i=1}^n |f^p_i + (\alpha_i g)^p|\Bigr)^{1/p}
\Bigr\|_X^p   
&=&
\sum_{i=1}^n \| |f^p_i +  (\alpha_i g)^p|^{1/p} \|_X^p   \nonumber \\
&\ge& (1- \varepsilon) \Bigl( \|g\|_X^p + \sum_{i=1}^n
\|f_i\|_X^p \Bigr). \label{eq3}
\end{eqnarray}

For (v) $\Rightarrow$ (vi) we apply the following separation
argument. Consider the convex set $B_{(X_{[p]})^*}$, which is a
compact Hausdorff space when endowed with the weak* topology, and
 the family of all functions $\Phi_{f_1,\dots ,f_n;
\alpha_1,\dots ,\alpha_n}\dopu  B_{(X_{[p]})^*} \to \mathbb{R}$, $n \in
\mathbb{N}$, $f_1,\dots ,f_n \in X_0$, $\alpha_1,\dots ,\alpha_n \in
\mathbb{R}$, $\sum_{i=1}^n \alpha^p_i=1$,  defined by
$$
\Phi_{f_1,\dots ,f_n; \alpha_1,\dots ,\alpha_n}(x^*):= (1- \varepsilon)
\Bigl( \|g\|^p + \sum_{i=1}^n \|f_i\|^p \Bigr) - \Bigl\langle
\sum_{i=1}^n |f_i^p + \alpha^p_i g^p|, x^* \Bigr\rangle .
$$

Each function defined in this way is clearly convex, and  the
family of all such functions is concave, since each convex
combination of two such functions can be written again as a
function of the same family; indeed for $0 \le \beta \le 1$
we have 
$$
\beta \Phi_{f_1,\dots ,f_n; \alpha_1,\dots ,\alpha_n} + (1-\beta)
\Phi_{\overline{f}_1,\dots ,\overline{f}_m;
\overline{\alpha}_1,\dots ,\overline{\alpha}_m}=
$$
$$
\Phi_{\beta^{1/p} f_1,\dots ,\beta^{1/p} f_n, (1-\beta)^{1/p}
  \overline{f}_1,\dots , 
(1-\beta)^{1/p} \overline{f}_m; \beta^{1/p} \alpha_1,\dots , \beta^{1/p} \alpha_n,
(1-\beta)^{1/p}\overline{\alpha}_1,\dots ,(1-\beta)^{1/p} \overline{\alpha}_m}
.
$$

These functions are continuous with respect to the weak* topology,
and by (\ref{eq3}) and the Hahn-Banach Theorem for each of them there is
an $x_1^* \in B_{(X_{[p]})^*}$ such that  $\Phi_{f_1,\dots ,f_n;
\alpha_1,\dots ,\alpha_n}(x_1^*) \le 0$, so an application of Ky
Fan's Lemma (see for instance \cite[Lemma~6.12]{libro}) gives an element
$x^*_0$ such that
$$
\Phi_{f_1,\dots ,f_n; \alpha_1,\dots ,\alpha_n}(x_0^*) \le 0
$$
for all the functions. Therefore, in particular, the inclusion in
(vi) is obtained.

Let us now prove (vi) $\Rightarrow$ (ii). Take a finite set of
elements $f_1,\dots ,f_n \in X$ and consider the finite dimensional
subspace $X_0$ generated by them. Take any slice
$S_{[p]}(x^*,\varepsilon_0)$ generated by a norm one element $x^*$
and an $\varepsilon >0$. Then an application of (vi) gives a $g
\in S^{1/p}_{[p]}(x^*,\varepsilon_0)$ and an element $x_0^* \in
B_{(X_{[p]})^*}$ such that
$$
\sum_{i=1}^n \langle |f^p_i+g^p|, x_0^* \rangle \ge (1-
\varepsilon)\Bigl( \sum_{i=1}^n\|f_i\|_X^p + n\|g\|_X^p\Bigr).
$$
Thus,
$$
\Bigl\| \Bigl(\sum_{i=1}^n |f_i|^p\Bigr)^{1/p}\Bigr\|_X^p + n \|g\|_X^p \ge (1-
\varepsilon)\Bigl( \sum_{i=1}^n\|f_i\|_X^p + n\|g\|^p_X\Bigr)
$$
and therefore,
$$
\Bigl\| \Bigl(\sum_{i=1}^n  |f_i|^p\Bigr)^{1/p}\Bigr\|^p_X + \varepsilon n \|g\|^p_X \ge (1-
\varepsilon)\Bigl( \sum_{i=1}^n\|f_i\|_X^p\Bigr).
$$
Since this construction can be done for every $\varepsilon >0$, we
obtain that
$$
\Bigl\| \Bigl(\sum_{i=1}^n  |f_i|^p\Bigr)^{1/p}\Bigr\|=
\Bigl(\sum_{i=1}^n\|f_i\|_X^p\Bigr)^{1/p}.
$$
Consequently, $X$ is an abstract $L^p$ space. Also,  for
$\varepsilon
>0$, taking a single function $f \in S_X$ and the subspace $X_0$
generated by it and an $x^* \in S_{(X_{[p]})^*}$ we obtain by (vi)
an $x^*_0 \in B_{(X_{[p]})^*}$ and a function $g \in
S^{1/p}_{[p]}(x^*,\varepsilon)$ such that
$$
 \||f^p + g^p|^{1/p} \|_X^p \ge \langle |f^p+g^p|, x_0^* \rangle 
\ge 2(1- \varepsilon) .
$$
Thus, Lemma~\ref{LemmapthDaugavet} gives that $X_{[p]}$ has the
Daugavet property.
\end{proof}

\begin{corollary} \label{corolfin}
Let $1 \le p < \infty$.
Every separable quasi-Banach function space satisfying the equivalent statements of Theorem~\ref{carac} is order isomorphic and isometric to $L^p([0,1])$.
\end{corollary}

This is a direct consequence of Theorem~\ref{carac} and the characterization of atomless separable $L^p$-spaces (see \cite[Theorem~2.7.3]{meyer-nieberg}).

\begin{remark}
Note that using Kakutani's representation theorem (see for instance
\cite[Theorem~1.b.2]{lint} or \cite[Theorem~2.7.1]{meyer-nieberg})
Theorem~\ref{carac} can be applied in a more abstract setting, without
the requirement for $X$ to be a quasi-Banach function space. If $X$
is just a Banach lattice 
that is also an abstract $L^p$ space, then
$X$  is order isometric to an $L^p(\mu)$ space over some measure
space $(\Omega,\Sigma,\mu)$, so in this case the condition of
$L^p(\mu)_{[p]} = L^1(\mu)$ having the Daugavet property, i.e., $\mu$
having no atoms, is characterized by the equivalent statemens of the
theorem. Therefore, Corollary~\ref{corolfin} can also be stated for
Banach lattices \textit{via} the atomic 
properties of the representing measure that Kakutani's theorem gives.
\end{remark}

\end{document}